\documentclass[leqno,12pt]{article}
\usepackage{amssymb,amsfonts}
\usepackage{amsmath,latexsym}

\newcommand{\bepr}{{\em Proof} } 
\newcommand{\enpr}{\hfill \rule{.5em}{.5em}}

\newcommand{\A}{{\mathbb A}}
 
\newcommand{\R}{{\mathbb R}}

\newcommand{\Tr}{{\rm Tr}}

\def\XXint#1#2#3{{\setbox0=\hbox{$#1{#2#3}{\int}$ }
\vcenter{\hbox{$#2#3$ }}\kern-.6\wd0}}

\newtheorem{defin}{Definition}[section] 
\newtheorem{prop}{Proposition}[section] 
\newtheorem{thm}{Theorem}[section] 
\newtheorem{lemma}{Lemma}[section]

\newtheorem{cor}{Corollary}[section] 
\newtheorem{rque}{Remark}[section] 
\newtheorem{notat}{Notation}[section]

\begin{document}

\title{Compensated Integrability in bounded domains~; Applications to gases}

\author{Denis Serre \\ \'Ecole Normale Sup\'erieure de Lyon\thanks{U.M.P.A., UMR CNRS--ENSL \# 5669. 46 all\'ee d'Italie, 69364 Lyon cedex 07. France. {\tt denis.serre@ens-lyon.fr}}}

\date{\today}

\maketitle

\begin{abstract}
An accurate functional inequality for Div-BV positive symmetric tensors $A$ in a bounded domain $U\subset\R^n$ arises whenever the tangential part of the normal trace $\gamma_\nu A\sim A\vec\nu$ is a finite measure over $\partial U$. The proof involves an extension operator to a neighbourhood of $\bar U$. The resulting inequality depends upon the domain only through the $C^3$-regularity of $\partial U$, some constant involving the curvature and its first derivatives.

This abstract statement applies to several models of Gas Dynamics (Euler system, Hard Spheres dynamics), as the boundary condition (slip, or reflection) tells us that $A\vec\nu$ is parallel to $\vec\nu$, where $A$ is the mass-momentum tensor.
\end{abstract}



\bigskip

\paragraph{Notations}

The domain $U\subset\R^n$ is open, bounded, with a $C^3$-boundary. The outer unit normal vector field to $\partial U$ is $\vec\nu\,:a\mapsto\vec\nu_a$. If $Z:\partial U\to\R^n$ is a vector field, its tangential component is the projection  $Z^\tau=Z-(Z\cdot\vec\nu)\vec\nu$. For a vector field $v:U\to\R^n$, the first differential is ${\rm d}v$ and the second one is ${\rm D}^2v$.

The space ${\cal M}(U)$ of finite Radon measures is understood in various senses, whether these measures are scalar, vector-valued or tensor-valued. The total mass (or {\em total variation} in some terminology) of a measure $\mu$ is $\|\mu\|_{\cal M}$. We consider also measures on the boundary $\partial U$. The bracket $\langle\mu,f\rangle$ stands for the duality, either between ${\cal M}(U)$ and $C_b(U)$, or between ${\cal D}'(U)$ and ${\cal D}(U)$.

For $n\times n$ matrices, $A:B=\Tr(B^TA)$ is the natural scalar product over ${\bf M}_n(\R)$. In practice, the matrices are symmetric, so the transposition doesn't matter.

An inequality ${\cal F}[U,A,\cdots]\le_n{\cal G}[U,A,\cdots]$ means that there exists a constant $c_n\in(0,+\infty)$, depending only upon the space dimension $n$, such that the inequality $${\cal F}[U,A,\cdots]\le c_n{\cal G}[U,A,\cdots]$$ holds true for every domain $U$, tensor $A$, ... under consideration.

Let $E$ be a normed vector space. If $u\in{\rm End}(E)$, we write $|u|_{\rm op}$ for the operator norm. Likewise, for $B\in{\rm Bil}(E;\R^n)$, we set
$$|B|_{\rm bil}=\sup_{x,y\in E}\frac{|B(x,y)|_{\R^n}}{|x|_E|y|_E}\,.$$
As usual $\|\cdot\|_\infty$ is a sup-norm.

%\vspace{3cm}

\section{Introduction and main results}

Given an open domain $U\subset\R^n$, we are concerned with symmetric $n\times n$ tensors $A$ whose entries, as well as the coordinates of the row-wise Divergence are finite measures:
$$\forall1\le i,j\le n,\quad a_{ij}\in{\cal M}(U),\quad\hbox{and}\quad\forall 1\le i\le n,\quad({\rm Div}\,A)_i:=\sum_{j=1}^n\partial_ja_{ij}\in{\cal M}(U).$$
These tensors form the Div-BV class. We make the important assumption that $A$ is positive semi-definite. In other words the measures 
$$A(\xi):=\sum_{i,j}\xi_i\xi_ja_{ij}$$
are non-negative, for every vector $\xi\in\R^n$.

When the domain $U$ is the whole space $\R^n$, it was shown  (see \cite{Ser_DPT,Ser_JMPA}) that the well-defined measure $(\det A)^{\frac1n}$, which is dominated by $\frac1n\,\Tr\,A$, is actually a measurable function of class $L^{\frac n{n-1}}$. This qualitative result (we speak of {\bf Compensated Integrability}), which says in particular that the singular part $A^{\rm sing}$ in the Radon-Nikodym decomposition satisfies $(\det A^{\rm sing})^{\frac1n}=0$ (this being also a consequence of an abstract result in \cite{ADHR}), is associated with a functional inequality
\begin{equation}\label{eq:FIRn}\tag{CI}
\int_{\R^n}(\det A)^{\frac1{n-1}}dx\le c_n\|{\rm Div}\,A\|_{{\cal M}(\R^n)}^{\frac n{n-1}}.
\end{equation}

We are interested here in the case where $U$ is bounded, with a smooth boundary. We ask whether the conclusion above still holds true for positive Div-BV tensors over $U$. The first observation is that if $\phi\in{\cal D}(U)$ is non-negative, then $\phi A$ is Div-BV, compactly supported and positive semi-definite, and thus its extension by $0_n$ to $U^c$ is Div-BV over $\R^n$. We infer that $(\det A)^{\frac1n}$ is still a measurable function, and belongs to $L^{\frac n{n-1}}_{\rm loc}(U)$. The estimate (\ref{eq:FIRn}), applied to $\widetilde{\phi A}$, provides an estimate for $(\det A)^{\frac1{n-1}}$ in $L^1_{\rm loc}(U)$, which deteriorates near the boundary $\partial U$.

An inequality such as (\ref{eq:FIRn}), with $U$ instead of $\R^n$,  is false because positive definite constant tensors violate it. 
To obtain an estimate of the full integral
$$\int_{U}(\det A)^{\frac1{n-1}},$$
we need some extra information about the normal trace $\gamma_\nu A$~; see Apprendix \ref{ap:trace} for a rigorous definition of this object, which mimics that in the space $H_{\rm div}(U)$ used in mathematical fluid dynamics \cite{Tem}. One naive strategy, already considered in \cite{Ser_DPT}, consists in extending $A$ by $0_n$ on $U^c$. The resulting tensor $\widetilde A$ satisfies, in the distributional sense,
$${\rm Div}\,\widetilde A=({\rm Div}\,A)|_U+\gamma_\nu A\,\left.{\cal H}_{n-1}\right|_{\partial U}.$$
In other words
$$\forall\vec\psi\in{\cal D}(\R^n;\R^n)\qquad \langle{\rm Div}\,\widetilde A,\vec\psi\rangle=\langle{\rm Div}\,A,\vec\psi|_U\rangle+\langle\gamma_\nu A,\vec\psi|_{\partial U}\rangle.$$
Since $\gamma_\nu A$ is only known to belong to the dual of $C^1(\partial U)$, $\widetilde A$ is not Div-BV in general, unless we assume explicitely that $\gamma_\nu A$ is a finite measure. In such a case, we have
$$\|{\rm Div}\,\widetilde A\|_{{\cal M}(\R^n)}=\|{\rm Div}\,A\|_{{\cal M}(U)}+\|\gamma_\nu A\|_{{\cal M}(\partial U)}.$$
This led us to state \cite{Ser_DPT}
\begin{prop}\label{p:faible}
Let $A\succ0_n$ be a Div-BV tensor over $U$, such that the normal trace $\gamma_\nu A$ is a (vector-valued) finite measure over $\partial U$. Then the measure $(\det A)^{\frac1n}$ actually belongs to $L^{\frac n{n-1}}(U)$, and there holds a functional inequality
\begin{equation}\label{eq:estfaible}
\int_U(\det A)^{\frac1{n-1}}dx\le_n\left(\|{\rm Div}\,A\|_{{\cal M}(U)}+\|\gamma_\nu A\|_{{\cal M}(\partial U)}\right)^{\frac n{n-1}}.
\end{equation}
\end{prop}
Alas the assumption  is too strong and the statement above could not be applied to inviscid gas dynamics. This contrasts with the case $U=\R^n$ for which (\ref{eq:FIRn}) led to a new fundamenetal {\em a priori} estimate~; see (\ref{eq:gasRd}) below.

\bigskip

An important remark is that the normal component of $\gamma_\nu A$, formally equal to $\vec\nu^TA\vec\nu$, must be a non-negative\footnote{A rigorous proof of this fact, not needed here, is a little bit technical.} distribution. As such, it is naturally a finite measure, so that the hypothesis upon the normal trace is nothing but the assumption that its tangential part 
$$\gamma_\nu^\tau A:=\gamma_\nu A-(\vec\nu\cdot\gamma_\nu A)\vec\nu$$ 
is a (vector-valued) finite measure. This suggests that a functional inequality such as (\ref{eq:estfaible}) could be improved by replacing the mass of $\gamma_\nu A$ by that of $\gamma_\nu^\tau A$. As we shall see below, this is almost true, in the sense that there is a reasonnable price to pay. This extra cost is of the form $K(\partial U)\|A\|_{\cal M}$ where $K(\partial U)^{-1}$ is some characteristic length of the domain.
Notice that this does not follow from an estimate of the mass of $\vec\nu^TA\vec\nu$ in terms of $\|{\rm Div}\,A\|_{\cal M}$ and $\|A\|_{\cal M}$, since the linear operator $A\mapsto\vec\nu^TA\vec\nu$ does not map the Div-BV space into ${\cal M}(\partial U)$ (that is, without the positivity assumption).
%Therefore the quantitative information
%about $\|\gamma_\nu^\tau A\|_{{\cal M}(\partial U)}$ is not sufficient to exploit (\ref{eq:estfaible})~; at this stage we still need the control of the full trace $\|\gamma_\nu A\|_{{\cal M}(\partial U)}$.

\bigskip

Improving the construction above, we show therefore that the knowledge of $\|\gamma_\nu^\tau A\|_{{\cal M}(\partial U)}$ essentially suffices to estimate the left-hand side of (\ref{eq:estfaible}), provided that the boundary $\partial U$ is of class $C^3$. The length scale mentionned above involves the following notations.
\begin{notat}\label{n:piu}
Let $\Sigma\subset \R^n$ be a closed $C^2$-hypersurface. We denote $L(\Sigma)$ the largest number $L$ such that every $x\in \R^n$ with $d(x;\Sigma)<L$ admits a unique projection on $\Sigma$. This projection is denoted $\pi_\Sigma$, or simply $\pi$.
\end{notat}
In the statements below, the hypersurface is the boundary $\partial U$ of our domain and is of class $C^3$. In most cases, $L(\partial U)$ is simply the infimum over  $\partial U$ of the curvature radii, equivalently $L(\partial U)=\|{\rm d}\vec\nu\|_\infty^{-1}$, but it might be smaller for dumbbell-shaped domains. For general domains, we only have $L^{-1}\ge\|{\rm d}\vec\nu\|_\infty$.
\begin{notat}
Let $\Sigma\subset \R^n$ be a closed $C^3$-hypersurface. We denote
$$K(\Sigma):=\max\left\{\frac1L\,,L\|{\rm D}^2\vec\nu\|_\infty\right\},\qquad L=L(\Sigma).$$
\end{notat}

\bigskip

Our main result reads as follows.
\begin{thm}\label{th:main}
Let $U\subset \R^n$ be a bounded open domain with $C^3$-boundary. 

For every symmetric positive semi-definite Div-BV tensor $A$ over $U$, whose tangential component $\gamma_\nu^\tau A$ of the normal trace belongs to ${\cal M}(\partial U;\R^n)$, the measure $(\det A)^{\frac1n}$  is actually a measurable function belonging to $L^{\frac n{n-1}}(U)$. The corresponding functional inequality is
\begin{equation}\label{eq:estfort}
\int_U(\det A)^{\frac1{n-1}}dx\le_n\left(\|{\rm Div}\,A\|_{{\cal M}(U)}+\|\gamma_\nu^\tau A\|_{{\cal M}(\partial U)}+K(\partial U)\|A\|_{{\cal M}(U)}\right)^{\frac n{n-1}}.
\end{equation}
\end{thm}

The proof of Theorem \ref{th:main} is based on the construction of a Div-BV extension of $A$ to the larger domain $U+B(0;L(\partial U))$. {\em A priori}, such an extension can be done for every Div-BV symmetric tensor, by following Babich strategy (see the seminal paper \cite{Bab}). Babich's idea provides an extended tensor whose Divergence does not charge $\partial U$. This is how trace and embedding theorems are proven in the context of Sobolev spaces. However a Babich-like extension does not preserve the positiveness of the tensor. Without this crucial condition, we may not apply Compensated Integrability to the extension of $A$. This is why we employ a more robust construction, which requires that $\gamma_\nu^\tau A$ is measure-valued.

\bigskip

For practical applications to evolution problems, we need the following variant of Theorem \ref{th:main} where the domain $U$ is of the form $\R\times\Omega$, the first factor standing for the time axis. The second factor, a physical domain, is a $d$-dimensional bounded open domain with $C^3$-boundary. We have $n=1+d$, so that the exponent $\frac1{n-1}$ is simply $\frac1d$\,. The independent variable is $x=(t,y)$.
\begin{thm}[$n=1+d$, $U=\R\times\Omega$.]\label{th:evol}
For every symmetric positive semi-definite Div-BV tensor $A$ over $\R\times\Omega$, whose tangential component $\gamma_\nu^\tau A$ of the normal trace belongs to ${\cal M}(\partial U;\R^n)$, the measure $(\det A)^{\frac1n}$  is actually a measurable function belonging to $L^{\frac n{n-1}}(U)$. The corresponding functional inequality is
\begin{equation}\label{eq:estevol}
\int_\R dt\int_\Omega(\det A)^{\frac1d}dy\le_n\left(\|{\rm Div}\,A\|_{{\cal M}(U)}+\|\gamma_\nu^\tau A\|_{{\cal M}(\partial U)}+K(\partial\Omega)\|A_{r}\|_{{\cal M}(U)}\right)^{1+\frac1d}.
\end{equation}
Hereabove, $A_{r}$ is the $n\times d$ right column  in the block-decomposition
\begin{equation}
\label{eq:blok}
A=\begin{pmatrix} A_{u\ell} & A_{ur} \\ A_{b\ell} & A_{br} \end{pmatrix}=(A_\ell\quad A_r).
\end{equation}
\end{thm}

Our applications are two-fold. On the one hand we consider an inviscid compressible gas obeying the Euler equations (conservation of mass and momentum, decay of mechanical energy). The mass-momentum tensor
$$F=\begin{pmatrix} \rho & \rho u^T \\ \rho u & \rho u\otimes u+pI_d \end{pmatrix},$$
where $\rho,p\ge0$ are the density, pressure and $u$ is the velocity field,
is positive and Div-free.  Because of the slip boundary condition $u\cdot\vec n=0$ on the boundary, the normal trace is 
$$\gamma_\nu F=\binom0{p\vec n}=p\vec\nu.$$
In particular $\gamma_\nu^\tau F\equiv0$ and Theorem \ref{th:evol} applies. With the help of (\ref{eq:estevol}), we establish the estimate
\begin{thm}\label{th:GD}
Assume a polytropic equation of state ($\varepsilon$ the specific internal energy, $\gamma>1$ the adiabatic constant)
\begin{equation}
\label{eq:poly}
p=(\gamma-1)\rho\varepsilon.
\end{equation}
Let an admissible inviscid gas flow in $\R_+\times\Omega$ have finite total mass $M$ and initial energy $E_0$ (the energy at time $t>0$ being $\le E_0$). Then we have for every $T>0$ 
\begin{equation}
\label{eq:estGD}
\int_0^Tdt\int_\Omega\rho^{\frac1d}p\,dy\le_d\sqrt{\frac{E_0}M}\,\left(M+\gamma T K(\partial\Omega)\sqrt{ME_0}\right)^{1+\frac1d}.
%\int_0^Tdt\int_\Omega\rho^{\frac1d}p\,dy\le_dM^{\frac1d}\left(\sqrt{ME_0\,}+\gamma K(\partial\Omega)TE_0\right).
\end{equation}
\end{thm}

\bigskip

\paragraph{Discussion.} The estimate (\ref{eq:estGD}) deteriorates when $T$ grows, since we expect the left-hand side to increase linearly with $T$, while the right-hand side goes superlinearly. It is more intersting when $T=T_f$ is the characteristic time of the flow, that is
\begin{equation}
\label{eq:charT}
T_f:=\frac1{\gamma K(\partial\Omega)}\,\sqrt{\frac M{E_0}}\,.
\end{equation}
Applying (\ref{eq:estGD}) with an arbitrary initial time $t_0>0$, we conclude
\begin{cor}
Under the same assumptions as in Theorem \ref{th:GD}, there holds
\begin{equation}
\label{eq:estGDTf}
\int_{t_0}^{t_0+T_f}dt\int_\Omega\rho^{\frac1d}p\,dy\le_dM^{1+\frac1d}\sqrt{\frac{E_0}M}\,,
\end{equation}
where the characteristic time of the flow is given by (\ref{eq:charT}).
\end{cor}

\bigskip

An interpretation of (\ref{eq:estGDTf}) is that in time-average, that is in Ces\`aro sense, we have
$$\int_\Omega\rho^{\frac1d}p\,dy\stackrel{\rm Ces.}= O\left(\gamma M^{\frac1d}K(\partial\Omega)E_0\right).$$
In other words,
$$\int_\Omega\rho^{\frac1d}p\,dy\stackrel{\rm Ces.}=O\left(\gamma \kappa(\partial\Omega)\bar\rho^{\frac1d}E_0\right),$$
where 
$$\bar\rho:=\frac M{{\rm vol}(\Omega)}\,,\qquad\kappa(\Omega):={\rm vol}(\Omega)^{\frac1d}K(\partial\Omega)$$ 
are the mean density and  a dimension-less expression which can be viewed as an {\em aspect ratio} of the domain, respectively. The latter is never small, because of  the isoperimetric inequality, but can be large for complicated domains.

The situation described above is significantly different from that of a gas filling the whole space $\R^d$, for which we obtained (see \cite{Ser_DPT}) the estimate
\begin{equation}
\label{eq:gasRd}
\int_0^{+\infty}dt\int_{\R^d}\rho^{\frac1d}p\,dy\le_dM^{\frac1d}\sqrt{ME_0\,}.
\end{equation}
The main difference between (\ref{eq:estGD}) and (\ref{eq:gasRd}) is that there was no privileged length or time scale for an inviscid\footnote{On the contrary, viscous gases always display characteristic length/time scales.} gas in absence of boundary. In particular the time integral in (\ref{eq:gasRd}) extends up to $T=+\infty$. This contrasts sharply with the case of a bounded domain~; say for instance that the gas is isentropic, so that (\ref{eq:poly}) reduces to $p=\alpha\rho^\gamma$ for some constant $\alpha>0$. Then by Jensen's Inequality,
$$\frac1{|\Omega|}\int_\Omega\rho^{\frac1d}p\,dy\ge\alpha\bar\rho^{\gamma+\frac1d}$$
is bounded away from zero as $t\to+\infty$, so that the time-integral in (\ref{eq:estGD}) diverges at $+\infty$.

\bigskip

Our second application concerns Hard Spheres Dynamics. The particles move with constant velocities between collisions, these involving either two particles or a particle against the boundary\footnote{It is costumary to discard the rare  configurations in which a collision involves three or more bodies.}. Collisions between particles are elastic, conserving the linear momentum and the kinetic energy. Collisions with the boundary are specular reflections. The change of velocity experienced by a ball at a collision is denoted $u_{\rm out}-u_{\rm in}$. We establish the following
\begin{thm}\label{th:HS}
Let $N$ identical balls of radius $r>0$ and mass $m>0$ move in $\Omega$, the collisions between particles being elastic, those between a particle and the boundary being specular reflections. Denote 
$$E=\frac m2\sum_{j=1}^N|u_j(t)|^2$$
the kinetic energy (a constant of the motion) and $\bar u=\sqrt{2E/Nm\,}$ the root mean square velocity. 

There exists a constant $\kappa_d>0$, depending neither upon the domain, nor on the initial configuration, such that if
\begin{equation}
\label{eq:seuil}
NrK(\partial\Omega)<\kappa_d,
\end{equation}
then we have for every $t\ge0$
\begin{equation}\label{eq:estHS}
\sum_{\rm coll.}^{(t,t+T_f)}|u_{\rm out}-u_{\rm in}|\le_dN^2\bar u.
\end{equation}
The sum above runs over all the collisions (either between particles, or against the boundary) occuring in the time interval $(t,t+T_f)$ where 
$$T_f:=\frac1{K(\partial\Omega)\bar u}$$
is the characteristic time of the flow.
\end{thm}

\paragraph{Remarks.}
\begin{itemize}
\item
The assumption and the estimate do not really depend upon the mass of the particles. Only their number and the root mean velocity matter. \item For reasonnably balanced domains, the threshold condition (\ref{eq:seuil}) means that the particles cannot fill more than some fraction of a diameter (or a bottleneck) of $\Omega$.  Even though this necessity is understandable, the probability that a large number of balls be approximately located on a single diameter is outstandingly small. Therefore we conjecture that the bound (\ref{eq:estHS}) remains valid for `most of' initial configurations under a much weaker assumption than (\ref{eq:seuil}). We leave as an {\bf Open Problem} to formulate and prove a quantitative version of this claim.
\item 
At reflections, $|u_{\rm out}-u_{\rm in}|$ is nothing but $2|u\cdot\vec n|$, so (\ref{eq:estHS}) provides a bound of 
$$\sum_{\rm refl.}^{(t_1,t_1+T_f)}|u\cdot\vec n|.$$
This does not seem good enough to evaluate the pressure at the boundary in the limit $N\to+\infty$, $Nm\to M$, since we expect
$$\int_{t_1}^{t_1+T_f}dt\int_{\partial\Omega} p\,dy=\lim\sum_{\rm refl.}^{(t_1,t_1+T_f)}m(u\cdot\vec n)^2.$$
At this point, we just lack a bound of the number of reflections. 
\end{itemize}

\bigskip

\paragraph{Plan of the article.}
Theorems \ref{th:main} and \ref{th:evol} are proved in Section \ref{s:abstract}. The proofs differ from each other only from one minor point, which is detailed in Paragraph \ref{sss:evol}.  Section \ref{s:GD} is dedicated to the application to the Euler system, with the proof of Theorem \ref{th:GD}. We treat Hard Spheres dynamics in Section \ref{s:HS}, where we recall in passing the notion of {\em determinantal masses} for tensors supported by graphs. Two appendices are dedicated to technical points. On the one hand, we recall the definition of the normal trace for Div-BV vector fields -- which applies to Div-BV tensors too. On the other hand, we establish classical estimates for the reflection across $\partial U$.

\bigskip

\paragraph{Acknowledgement.} This research was carried out, for a part, during a stay of the author at the Korean Advanced Institute of Science \& Technology, Daejon at the kind invitation of Prof. Moon-Jin Kang.

\section{Push-forward and extension of Div-BV tensors}\label{s:abstract}

\subsection{Push-forward by a diffeomorphism}

Let $\omega, V\subset\R^n$ be open domains, and $\phi:\omega\to V$ be a $C^2$-diffeomorphism. We use latin indices for the coordinates in $V$ and greek ones in $\omega$. 

If $A:\omega\to{\bf Sym}_n$ is a symmetric tensor, whose entries are Radon measures, we define a push-forward symmetric tensor $\phi_*A$ over $V$ by the formula
$$\langle \phi_*A,M\rangle_V:=\langle{\rm d}\phi^TA\,{\rm d}\phi,M\circ\phi\rangle_\omega$$
for every $M\in C_c(V;{\bf M}_n(\R))$. For completeness,  the $(i,j)$-entry of ${\rm d}\phi^TA\,{\rm d}\phi$ is
$$\sum_{\alpha,\beta}a_{\alpha\beta}\partial_\alpha\phi_i\partial_\beta\phi_j,$$
so that the map $A\mapsto\phi_*A$ preserves positive semi-definiteness.
\begin{lemma}
The Divergence of $\phi_*A$ is given by
$$\langle {\rm Div}\,(\phi_*A),Y\rangle_V
=\langle({\rm Div}\,A){\rm d}\phi+A:{\rm D}^2\phi,Y\circ\phi\rangle_\omega,$$
for every $Y\in{\cal D}(V;\R^n)$.
\end{lemma}
As above, the $i$-th coordinate of $({\rm Div}\,A){\rm d}\phi+A:{\rm D}^2\phi$ is nothing but $({\rm Div}\,A)\cdot\nabla\phi_i+\Tr(A{\rm D}^2\phi_i)$, that is ${\rm div}(A\nabla\phi_i)$. The Lemma implies
\begin{cor}
If moreover $d\phi$ is globally Lipschitz, then 
$$(A\in{\rm DivBV}(\omega))\Longrightarrow(\phi_*A\in{\rm DivBV}(V)),$$ 
with an estimate
\begin{equation}
\label{eq:pfDiv}
\|{\rm Div}(\phi_*A)\|_{{\cal M}(V)}\le\|d\phi\|_\infty\|{\rm Div}\,A\|_{{\cal M}(\omega)}+\|{\rm D}^2\phi\|_\infty\|A\|_{{\cal M}(\omega)}.
\end{equation}
\end{cor}

Concerning the trace, we have
\begin{prop}
Let $\phi$ be a $C^2$-diffeomorphism such that $d\phi$ is globally Lipschitz, and let $A$ be Div-BV over $\omega$. Then the normal trace of $\phi_*A$ is given by the formula
\begin{equation}\label{eq:trphiA}
\langle\gamma_\nu(\phi_*A),Z\rangle_{\partial V}=\langle\gamma_\nu(Ad\phi),Z\circ\phi\rangle_{\partial\omega}.
\end{equation}
\end{prop}

The right-hand side of (\ref{eq:trphiA}) can be developed as
$$\sum_{i=1}^n\langle\gamma_\nu(A\nabla\phi_i),Z_i\circ\phi\rangle_{\partial\omega},$$
where each vector field $A\nabla\phi_i$ is div-BV. Notice that because of (\ref{eq:Cungam}), we have
\begin{equation}
\label{eq:Trefle}
\gamma_\nu(A\nabla\phi_i)=(\gamma_\nu A)\cdot\nabla\phi_i.
\end{equation}

\bigskip

\bepr

We proceed by duality. For compactly supported test vector field $Y$, the chain rule gives
\begin{eqnarray*}
\langle {\rm Div}\,(\phi_*A),Y\rangle_V & = & -\langle \phi_*A,\nabla Y\rangle_V=-\langle A\,{\rm d}\phi,\nabla(Y\circ\phi)\rangle_\omega \\
& = & \langle{\rm Div}(A\,{\rm d}\phi),Y\circ\phi\rangle_\omega
=\langle({\rm Div}\,A){\rm d}\phi+A:{\rm D}^2\phi,Y\circ\phi\rangle_\omega.
\end{eqnarray*}
If moreover $|Y(x)|\le1$ everywhere, then $|Y\circ\phi|\le1$ as well, and we have
$$|\langle {\rm Div}\,(\phi_*A),Y\rangle_V|
\le\|({\rm Div}\,A){\rm d}\phi+A:{\rm D}^2\phi\|_{{\cal M}(\omega)}
\le\|d\phi\|_\infty\|{\rm Div}\,A\|_{{\cal M}(\omega)}+\|{\rm D}^2\phi\|_\infty\|A\|_{{\cal M}(\omega)},$$
whence (\ref{eq:pfDiv}).

Finally, if $Y$ is not  compactly supported,
\begin{eqnarray*}
\langle\gamma_\nu(\phi_*A),Y|_{\partial V}\rangle & = & \langle\phi_*A,\nabla Y\rangle_V+\langle {\rm Div}\,(\phi_*A),Y\rangle_V \\
& = & \langle A\,{\rm d}\phi,\nabla(Y\circ\phi)\rangle_\omega + \langle{\rm Div}(A\,{\rm d}\phi),Y\circ\phi\rangle_\omega \\
& = & \langle\gamma_\nu(Ad\phi),(Y\circ\phi)|_{\partial\omega}\rangle=\langle\gamma_\nu(Ad\phi),Y|_{\partial V}\circ\phi\rangle.
\end{eqnarray*}

\enpr

\subsection{Extension across a boundary}\label{ss:ext}

Let $\Sigma$ be a closed $C^3$-hypersurface  in $\R^n$, with unit normal vector field $a\to\vec\nu_a$. The open set $\Sigma_L\subset\R^n$ defined by $d(x;\Sigma)<L(\Sigma)$ is split by $\Sigma$ in two halves denoted $\omega$ and $V$. The smooth involution $\sigma:\Sigma_L\to\Sigma_L$ given by
$$\sigma(x)=2\pi_\Sigma(x)-x$$
exchanges $\omega$ and $V$, while fixing the points of $\Sigma$. For each point $a\in\Sigma$, the segment $(a-L\vec\nu_a,a+L\vec\nu_a)$ is stable under $\sigma$, which acts as $a+s\vec\nu_a\mapsto a-s\vec\nu_a$.

Let $A$ be a symmetric Div-BV tensor over $\omega$. We have at our disposal a Div-BV tensor over $V$, namely $\sigma_*A$. This allows us to define an extension $A_{\rm ext}$ to $\Sigma_L$~: its restrictions to $\omega$ and $V$ are $A$ and $\sigma_*A$, respectively, and it does not charge $\Sigma$. In other words
$$\forall M\in C_c(\Sigma_L;{\bf M}_n(\R)),\qquad\langle A_{\rm ext},M\rangle
=\langle A,M|_\omega\rangle+\langle \sigma_*A,M|_V\rangle.$$
Notice that if $A$ is positive semi-definite, then so is $A_{\rm ext}$.
The natural question of whether $A_{\rm ext}$ is still Div-BV, reduces to: how much ${\rm Div}\,A_{\rm ext}$ charges $\Sigma$~?
\begin{prop}\label{p:extBV}
If the symmetric tensor $A$ is Div-BV over $\omega$, and if the tangential component $\gamma_\nu^\tau A$ of is normal trace is a (vector-valued) finite measure over $\Sigma$, then $A_{\rm ext}$ is Div-BV over $\Sigma_L$, and we have
\begin{equation}
\label{eq:Divext}
\langle {\rm Div}\,A_{\rm ext},\vec\psi\rangle = \langle {\rm Div}\,A,\vec\psi|_\omega\rangle+\langle {\rm Div}\,(\sigma_*A),\vec\psi|_V\rangle-2\langle\gamma_\nu^\tau A,\vec\psi|_\Sigma\rangle.
\end{equation}
\end{prop}

\bigskip

\bepr

For definiteness, we assume that the normal vector field $\vec\nu$ points towards $V$.

That $A_{\rm ext}$ is a (matrix-valued) finite measure is clear. We only have to evaluate its Divergence. For test vector fields, we have
\begin{eqnarray*}
\langle {\rm Div}\,A_{\rm ext},\vec\psi\rangle & = & -\langle A_{\rm ext},\nabla\vec\psi\rangle=-\langle A,\nabla\vec\psi|_\omega\rangle-\langle\sigma_*A,\nabla\vec\psi|_V\rangle \\
& = & \langle {\rm Div}\,A,\vec\psi|_\omega\rangle-\langle \gamma_\nu A,\vec\psi|_\Sigma\rangle+\langle {\rm Div}\,(\sigma_*A),\vec\psi|_V\rangle-\langle \gamma_\nu(\sigma_*A),\vec\psi|_\Sigma\rangle
\end{eqnarray*}
Involving (\ref{eq:trphiA}) and (\ref{eq:Cungam}), and using the fact that $\sigma|_\Sigma$ is the identity, together with (\ref{eq:Trefle}), we have
\begin{eqnarray*}
\langle \gamma_\nu A,\vec\psi|_\Sigma\rangle+\langle \gamma_\nu(\sigma_*A),\vec\psi|_\Sigma\rangle & = & \langle \gamma_\nu(A(I_n+{\rm d}\sigma)),\vec\psi|_\Sigma\rangle \\
& =  & \langle (I_n+{\rm d}\sigma)\gamma_\nu A,\vec\psi|_\Sigma\rangle.
\end{eqnarray*}
There remains to notice that the restriction of $d\sigma$ to $\Sigma$ is the orthogonal symmetry $I_n-2\vec\nu\otimes\vec\nu$, and thus $I_n+d\sigma_a$ is twice the orthogonal projection on $\vec\nu_a^\bot$, the tangent space. Hence
$$\langle \gamma_\nu A,\vec\psi|_\Sigma\rangle+\langle \gamma_\nu(\sigma_*A),\vec\psi|_\Sigma\rangle=2\langle\gamma_\nu^\tau A,\vec\psi|_\Sigma\rangle.$$

\enpr

\begin{rque}
We could have used another diffeomorphism $\phi$ instead of $\sigma$, with the property that $\phi$ exchanges $\omega$ and $V$, and it fixes $\Sigma$. However, the requirement that $d\phi=I_n-2\vec\nu\otimes\vec\nu$ along $\Sigma$ is a rather strong one. It rules out simpler $C^2$-flips associated with less regular manifolds. An important example is the ``horizontal flip'', a kind of nonlinear transvection: say that $\Sigma$ is locally given by an equation $x_n=f(y)$ where $y=(x_1,\ldots,x_{n-1})$ and $f\in C^2$. The horizontal flip is the map
$$\phi\binom y{x_n}=\binom y{2f(y)-x_n}\,.$$
Then at $a=(y,f(y))$,
$$I_n+d\phi(a)=2\begin{pmatrix} I_{n-1} & 0 \\ df(y) & 0 \end{pmatrix}$$
and the expression $(I_n+{\rm d}\sigma)\gamma_\nu A$ involves the horizontal component (that along the coordinate $y$) of $\gamma_\nu A$. This component, being frame-dependent, is useless in applications such as gas dynamics, where it has no physical meaning.
\end{rque}

\subsection{Proof of Theorems \protect\ref{th:main} and \protect\ref{th:evol}}

Let $A\succ0_n$ be Div-BV in the bounded domain $U\subset\R^n$, such that $\gamma_\nu^\tau A$ is measure-valued. We apply the construction above with $\Sigma=\partial U$, and we write $L$ for $L(\partial U))$. The domain
$$U\cup\Sigma_L= U+B(0;L),$$
is denoted $U_L$.

One of $\omega$ and $V$, say $\omega$, is included in $U$, while the other is exterior to $U$. The restriction $A|_\omega$ admits a Div-BV extension $(A|_\omega)_{\rm ext}$ to $\Sigma_L$. Since it matches $A$ in $\omega$, we have a Div-BV extension of $A$ to $U_L$,
$$A'=\left\{\begin{array}{lcr}
A & \hbox{in} & U, \\
(A|_\omega)_{\rm ext} & \hbox{in} & \Sigma_L.
\end{array}\right.$$
Differential operators being local, we infer from (\ref{eq:Divext}) the formula
$$\langle {\rm Div}\,A',\vec\psi\rangle = \langle {\rm Div}\,A,\vec\psi|_U\rangle+\langle {\rm Div}\,(\sigma_*A|_\omega),\vec\psi|_V\rangle-2\langle\gamma_\nu^\tau A,\vec\psi|_\Sigma\rangle,$$
from which there comes
\begin{equation}
\label{eq:DivApr}
\|{\rm Div}\,A' \|_{\cal M}=\|{\rm Div}\,A \|_{\cal M}+\|({\rm Div}\,(\sigma_*A|_\omega))|_V \|_{\cal M}+\|\gamma_\nu^\tau A \|_{\cal M}.
\end{equation}

Let now $\chi\in{\cal D}(U_L)$ be such that $0\le\chi\le1$ in $U_L$,  $\chi\equiv1$ over $U$, and $|\nabla\chi|<\frac4L$ otherwise. The tensor $\chi A'$ is still Div-BV, positive semi-definite, and in addition compactly supported in $U_L$. Thus its extension $\widetilde{\chi A'}$ by $0_n$ is Div-BV in $\R^n$. We have
\begin{eqnarray*}
\|{\rm Div}\,\widetilde{\chi A'} \|_{\cal M} & = & \|{\rm Div}\,\chi A' \|_{\cal M}\le \|\chi\|_\infty\|{\rm Div}\,A' \|_{\cal M}+\|\nabla\chi\|_\infty\|A' \|_{\cal M} \\
& \le & \|{\rm Div}\,A \|_{\cal M}+\|({\rm Div}\,(\sigma_*A|_\omega))|_V \|_{\cal M}+\|\gamma_\nu^\tau A \|_{\cal M}+\frac4L\,\|A' \|_{\cal M}  \\
& \le & (1+\|{\rm D}\sigma\|_\infty)\|{\rm Div}\,A \|_{\cal M}+\|\gamma_\nu^\tau A \|_{\cal M}+\left(\frac4L+\|{\rm D}^2\sigma\|_\infty\right)\,\|A' \|_{\cal M},
\end{eqnarray*}
where we used (\ref{eq:pfDiv}).

To make this estimate fully explicit, we need bounds of $\|{\rm D}\sigma\|_\infty$ and $\|{\rm D}^2\sigma\|_\infty$. For this, we make a careful choice for the cut-off $\chi$~: it vanishes whenever $d(x;\partial U)>\frac L2$\,, this being compatible with the bound $|\nabla\chi|<\frac4L$\,. Therefore the $L^\infty$-bound of derivatives of the involution $\sigma$ need only be taken over the smaller domain
$$\omega_{\frac12}=\left\{x\in U\,|\,0<d(x;\partial U)<\frac L2\right\}.$$
These are evaluated in Appendix \ref{ap:refl}, where we obtain (see Lemma \ref{lm:hddh} and Corollary \ref{c:ddsig})
$$\|{\rm d}\sigma\|_{L^\infty(\omega_{1/2})}\le 3,\qquad \|{\rm D}^2\sigma\|_{L^\infty(\omega_{1/2})}\le 16\,K(\partial U).$$
We conclude therefore
\begin{equation}
\label{eq:estext}
\|{\rm Div}\,\widetilde{\chi A'} \|_{\cal M}  \le  4\,\|{\rm Div}\,A \|_{\cal M}+\|\gamma_\nu^\tau A \|_{\cal M}+20\,K(\partial U)\,\|A' \|_{\cal M}.
\end{equation}

We are now in position to apply the functional inequality (\ref{eq:FIRn}) to $\widetilde{\chi A'}$. Since $\det\widetilde{\chi A'}$ is non-negative and its restriction to $U$ equals $\det A$, we obtain
\begin{eqnarray*}
\int_U(\det A)^{\frac1{n-1}}dx & \le & \int_{\R^n}(\det\widetilde{\chi A'})^{\frac1{n-1}}dx\le_n\|{\rm Div}\,\widetilde{\chi A'}\|_{\cal M}^{\frac n{n-1}} \\
& \le_n & \left(\,\|{\rm Div}\,A \|_{\cal M}+\|\gamma_\nu^\tau A \|_{\cal M}+K(\partial U)\,\|A' \|_{\cal M}\right)^{\frac n{n-1}}.
\end{eqnarray*}
With 
$$\|A' \|_{\cal M}\le\|A \|_{\cal M}+\|\sigma_*A|_{\omega_{1/2}} \|_{\cal M}\le\|A \|_{\cal M}+\|{\rm d}\sigma\|_\infty^2\|A \|_{\cal M}\le10\|A \|_{\cal M},$$
this ends the proof of Theorem \ref{th:main}. 

\subsubsection{The case  of Theorem \protect\ref{th:evol}}\label{sss:evol}

When $U=\R\times\Omega$ (and $x=(t,y)$) instead, we proceed the same way, with the following adaptations. The hypersurface $\Sigma=\R\times\partial\Omega$ splits $U_L$ (remark the equality $L(\partial U)=L(\partial\Omega)$) into $\R\times\omega$ and $\R\times V$ where $\omega\subset\Omega$ and $V$ is exterior to $\Omega$. The involution is $\sigma={\rm id}_\R\otimes\sigma_\omega$ where $\sigma_\omega:\omega\to V$ is defined as above. The only new fact is the evaluation of the contribution of $\|A\|_{\cal M}$ in the functional inequality. We recall that it originates from the terms $A:{\rm D}^2\sigma$ and $\tilde A\nabla\chi$.

On the one hand we have
we have
$$A:{\rm D}^2=A_{u\ell}\partial_t^2+2(A_{b\ell}\cdot\nabla_y)\partial_t+A_{br}:{\rm D}_{yy}^2,$$
where we use the block decomposition (\ref{eq:blok}).
Since $\partial_t^2\sigma\equiv0$ and $\nabla_y\partial_t\sigma\equiv0$, this yields
$$A:{\rm D}^2\sigma=\binom0{A_{br}:{\rm D}_y^2\sigma_\omega},$$
and therefore
$$\|A:{\rm D}^2\sigma\|_{\cal M}=\|A_{br}:{\rm D}_y^2\sigma_\omega\|_{\cal M}\le\|{\rm D}_y^2\sigma_\omega\|_\infty\|A_{br}\|_{\cal M}.$$

On the other hand the cut-off function $\chi$ depends only upon the space variable $y$, and satisfies as above $|\nabla\chi|\le\frac4L$\,. It contributes to ${\rm Div}\,\widetilde{\chi A'}$ through the term $A'\nabla\chi$, which involves only the $n\times d$ block $(A')_r$. Whence a contribution (at most)
$$\frac4L\,\|(A')_r\|_{\cal M}$$
to the mass of the Divergence. However, noticing that the differential $d\sigma={\rm diag}(1,d\sigma_\omega)$ is block-diagonal, we see that $(A')_r$ depends only upon the right block $A_r$, and not upon the left one (that is, not upon the top-left entry $a_{00}$).  The contribution above is thus bounded by 
$$\frac{16}L\,\|A_r\|_{\cal M}.$$

This explains why the factor $\|A\|_{\cal M}$ in (\ref{eq:estfort}) is replaced by $\|A_r\|_{\cal M}$ in (\ref{eq:estevol}).

\section{Application to inviscid gas dynamics}\label{s:GD}

In $d$-dimensional gas dynamics, the tensor
$$F=\begin{pmatrix} \rho & \rho u^T \\ \rho u & \rho u\otimes u+pI_d \end{pmatrix}$$
is symmetric, positive (because $p\ge0$) and Div-free (because of conservation of mass and momentum). 
The boundary $\partial\Omega$ being impermeable, the boundary condition is 
$$u\cdot\vec n\equiv0\qquad\hbox{over }\Gamma=(0,T)\times\partial\Omega.$$

Let $(0,T)$ be a time interval. We define a tensor $A$ over $U=\R\times\Omega$ by
$$A=\left\{\begin{array}{lcr}
F & \hbox{if} & t\in(0,T), \\
0_n & \hbox{if not.} & 
\end{array}\right.$$
Denoting $\vec n$ the outer unit normal vector to $\partial\Omega$,  the normal trace $\gamma_\nu A$ over $\Gamma$ is
$$\binom0{p\vec n},$$
so that $\gamma_\nu^\tau A\equiv0.$
Therefore Theorem \ref{th:evol} applies, and the contribution of the boundary vanishes. With $\det F=\rho^dp$, we have 
$$\int_0^Tdt\int_\Omega \rho^{\frac1d}p\,dy\le_n\left(\|{\rm Div}\,A\|_{\cal M}+K(\partial\Omega)\|A_{r}\|_{\cal M}\right)^{1+\frac1d}.$$
Since $F$ is Div-free, the Divergence of $A$ reduces to
$${\rm Div}\,A=\binom{\rho(T)}{\rho u(T)}{\cal H}_d|_{t=T}-\binom{\rho(0)}{\rho u(0)}{\cal H}_d|_{t=0},$$
where ${\cal H}_d$ stands for the $d$-dimensional Hausdorf measure, here the Lebesgue measure over hyperplanes. With Cauchy-Schwarz Inequality, we infer
$$\|{\rm Div}\,A\|_{\cal M}\le 2M+\sqrt{2ME_0\,}+\sqrt{2ME(T)\,}\le 2(M+\sqrt{2ME_0\,}),$$
where $M=\int_\Omega\rho(t,y)\,dy\equiv\int_\Omega\rho_0(y)\,dy$ is the mass of the gas, and 
$$E_0=\int_\Omega\left(\frac12\rho|u|^2+\rho\varepsilon\right)(0,y)\,dy$$
is the mechanical energy at initial time.

On another hand the formul\ae\, $F_{tr}=\rho u^T$ and $F_{br}=\rho u\otimes u+pI_d$ yield
$$\|A_{r}\|_{\cal M}\le\int_0^Tdt\int_\Omega(\rho |u|+\rho|u|^2+dp)\,dy.$$
Assuming the polytropic equation of state (\ref{eq:poly}), this gives
$$\|A_{r}\|_{\cal M}\le T\left(\sqrt{2ME_0}+\max\{2,d(\gamma-1)\}E_0\right).$$
We conclude that an admissible flow satisfies the estimate
\begin{equation}
\label{eq:nonhom}
\int_0^Tdt\int_\Omega \rho^{\frac1d}p\,dy\le_n\left(M+\sqrt{ME_0\,}+TK(\partial\Omega)\left(\sqrt{ME_0}+\gamma E_0\right)\right)^{1+\frac1d}.
\end{equation}
This one has a flaw however, because it is not homogeneous from the point of view of physical dimensions. For instance $M$ and $\sqrt{ME_0\,}$ have different dimensions. To achieve the homogeneous estimate (\ref{eq:estGD}), we observe that for every parameter $\mu>0$, the change of dependent/independent variables
$$(t,y,\rho,u,p,\varepsilon)\longmapsto(t',y',\rho',u',p',\varepsilon')=(\frac t\mu\,,y,\rho,\mu u,\mu^2p,\mu^2\varepsilon)$$
defines another flow, with same domain $\Omega$, same mass $M'=M$, but energy $E_0'=\mu^2E_0$. Applying (\ref{eq:nonhom}) to this flow, on the time interval $(0,\frac T\mu)$, we obtain
$$\int_0^{T/\mu}dt'\int_\Omega\rho^{\frac1d}(\mu^2 p)\,dy\le_n\left(M+\mu\sqrt{ME_0\,}+T K(\partial\Omega)\left(\sqrt{ME_0}+\mu \gamma E_0\right)\right)^{1+\frac1d}.$$
Returning to the time variable $t$ in the integral, this gives us a family of inequalities, parametrized by $\mu\in(0,+\infty)$,
$$\mu\int_0^Tdt\int_\Omega \rho^{\frac1d}p\,dy\le_n\left(M+\mu\sqrt{ME_0\,}+T K(\partial\Omega)\left(\sqrt{ME_0}+\mu \gamma E_0\right)\right)^{1+\frac1d}.$$
We now select the parameter $\mu=(M/E_0)^{\frac12}$, the inverse of a characteristic velocity, to conclude
$$\int_0^Tdt\int_\Omega \rho^{\frac1d}p\,dy\le_n\frac1\mu\,\left(M+\gamma T K(\partial\Omega)\sqrt{ME_0}\right)^{1+\frac1d}=\sqrt{\frac{E_0}M}\,\left(M+\gamma T K(\partial\Omega)\sqrt{ME_0}\right)^{1+\frac1d}.$$

This ends the proof of Theorem \ref{th:GD}.

\section{Application to Hard Spheres dynamics}\label{s:HS}

In Hard Spheres dynamics, the mass-momentum tensor $M$ (see \cite{Ser_HS}) is singular, supported by a graph $\cal G$ made of the trajectories of the centers of particles,  completed by the horizontal segments between the centers of colliding particles, at collision times. Let us recall that a {\em colliton}, a contribution to $M$, is the tensor accounting for the exchange of linear momentum between two colliding particles. Denoting $\tau\in(0,T)$ the instant of the collision, and $p,p'\in\Omega$ the positions of the colliding particles (one has $|p'-p|=2r$), it takes the form
$$\frac m{|[v]|}\,\binom0{[v]}\otimes\binom0{[v]}\,d\ell_J,\qquad J:=\{\tau\}\times(p,p'),$$
where $[v]=v_2-v_1$ is the difference between the ante/post velocities of one of both particles, and $d\ell_J$ is the length element on the interval $J$. Notice that because of $v_2'-v_1'=v_1-v_2$ (conservation of momentum), it does not matter which of the colliding particles is used when constructing the colliton. We recall in passing that $[v]$ is parallel to $(p,p')$.

The Div-freeness of $M$ expresses the conservation of particles and of linear momentum. This tensor being esentially rank-one, its determinant vanishes and thus (\ref{eq:FIRn}) is useless in this context. One must add correctors at each node of $\cal G$, and then apply an appropriate version of Compensated Integrability, which incorporates the so-called {\em determinantal masses}, a notion introduced in \cite{Ser_HS}. Needless to say, this version is valid in the context of bounded domains, the proof of the augmented functional inequality being exactly the same as in our recent work \cite{Ser_period}. 

We make generic assumptions, under which the dynamics is globally defined:
\begin{itemize}
\item Collisions happen either between pairs of particles, or between a single particle and the boundary, 
\item Collisions form a discrete set (they do not accumulate),
\item The radius $r$ of the balls is smaller than $L(\partial\Omega)$, so that a particle colliding the boundary has only one contact point with $\partial\Omega$. Indeed we shall need a stronger assumption~; see condition (\ref{eq:seuil}).
\end{itemize}

Since we are interested in motions in the bounded domain $\Omega$, surrounded by a rigid boundary, we extend the notion of {\em colliton} to the case of a single particle bouncing against the wall. This will ensure that $M$ is a Div-free tensor. Such collitons are defined as above, but now $p'$ is replaced by the contact point $q\in\partial\Omega$ of the particle. Because of the specular law of reflection, we know that $[v]$ is colinear to the normal $\vec n_q$, hence to $\vec{pq}$. Since the unit normal to $\Gamma$ at $(\tau,q)$ is $\vec\nu=\binom0{\vec n_q}$, the colliton may be recast as
$$m|[v]|\vec\nu\otimes\vec\nu\,d\ell_J,\qquad J:=\{\tau\}\times(p,q),$$
and its contribution to the normal trace $\gamma_\nu M$ is a Dirac mass
\begin{equation}
\label{eq:bdrycoll}
m\,\binom0{[v]}\,\delta_{(\tau,q)}=m|[v]|\vec\nu\,\delta_{(\tau,q)}.
\end{equation}

Since the nodes of $\cal G$ stay away from $\Gamma:=(0,T)\times\partial\Omega$ (by a distance $r$), and the supports of the correctors can be chosen arbitrarily small, thus disjoint from $\Gamma$, the normal trace of these correctors vanishes identically. Thus only the mass-momentum tensor matters along the lateral boundary.  The situation is actually the same for the trajectories of the centers of mass, which stay away from $\Gamma$. Thus only the boundary collitons matter when computing the normal trace $\gamma_\nu M$ along $\Gamma$. From formula (\ref{eq:bdrycoll}), this tangential part vanishes:
$$\gamma_\nu^\tau M\equiv0.$$

\paragraph{Remarks.} Our proof of Theorem \ref{th:HS} will use the  extension of $M$ as a Div-BV tensor over $U_L$, as constructed in Section \ref{ss:ext}. At least two other extensions can be constructed easily. But although being  more explicit than the one designed in Section \ref{ss:ext}, they do not seem to fit well with Compensated Integrability.
\begin{itemize}
\item A first option consist in extending the support of a boundary colliton beyond the contact point $q$. The resulting extension is Div-free. But applying C.I. with this special extension yields nowhere. Mimicking the calculations below, we end up with an inequality of the form
$$\sum_{\rm coll.}|u_{\rm out}-u_{\rm in}|\le_d N\left(N\bar u+\sum_{\rm coll.}|u_{\rm out}-u_{\rm in}|\right),$$
which does not imply an estimate of $\sum_{\rm coll.}|u_{\rm out}-u_{\rm in}|$.
\item An other construction follows from the remark that a particle $\pi$ bouncing against the boundary behaves as if a virtual, point particle $\pi'$ of same mass, coming from outside at velocity $u_{\rm out}$, collides with $\pi$ at $(\tau,q)$ (with the same notations as above). Then $\pi'$ bounces back with outgoing velocity $u_{\rm in}$ ($\pi$ and $\pi'$ exchange their velocities in  the collision). Then the extension $\widetilde M$, still Div-free, is the mass-momentum tensor accounting for both real and virtual particles. Equivalently, $N$ is replaced by the number $N'$ or real and virtual particles, a quantity we are unable to control.
\end{itemize}

\bigskip

It is therefore better to deal with an extension of $M$ as constructed in Section \ref{ss:ext} for general Div-BV symmetric tensors. 

\subsection{Determinantal masses}

The concept of {\em determinantal masses} was introduced in \cite{Ser_HS}. However we refer to our recent paper \cite{Ser_period} for a more elegant presentation.

Let us recall that if $E\subset\R^n$ is a linear subspace of dimension $\ell$, then the $\ell$-dimensional Lebesgue measure ${\cal L}_\ell|_E$ is a homogeneous distribution, of order $\ell-n$. For instance, the Dirac mass $\delta_0$ is homogeneous of degree $-n$. We shall also speak of homogeneous distributions about a point $X$, when homogeneity is restored after the translation $x\mapsto x-X$.

Let $A\succ0_n$ be a Div-BV tensor in $\R^n$, which is Div-free and positively homogeneous of degree $1-n$ in a neighbourhood of $X\in \R^n$. It was showed in \cite{Ser_JMPA} that the latter property is somewhat extreme in the realm of Div-free tensors: there exists a non-negative measure $\mu$ on the unit sphere such that 
$$A=\mu\left(\frac{x-X}r\right)\,\frac{(x-X)\otimes(x-X)}{r^{n+1}}\,,\qquad r=|x-X|,$$
and
$$\int_{S_{n-1}}\vec e d\mu(\vec e)=0.$$
Then Pogorelov's Theorem \cite{Pog} about the Minkowski Problem ensures that there exists a convex function $\theta$, positively homogeneous (about $X$) of degree $1$, such that
$$A=\widehat{{\rm D}^2\theta}$$
in this neighbourhood. We recall the
\begin{defin}
The {\em determinantal mass} of $A$ at $X$, denoted ${\rm Dm}(A;X)$, is the volume of the convex body $\partial\theta(X)$ (the subgradient).
\end{defin}
So far, the expression $(\det A)^{\frac1{n-1}}$, viewed as the $\frac n{n-1}$'th power of $(\det A)^{\frac1n}$ was proved to be an integrable function, see (\ref{eq:FIRn}). The variant eleborated in \cite{Ser_HS} suggests that this function is only the absolutely continuous part of a non-negative measure, whose singular part contains a Dirac mass at $X$, of weight ${\rm Dm}(A;X)$.

We recall now how determinantal masses incorporate into Compensated Integrability.
\begin{thm}[D. S. \cite{Ser_HS}]\label{th:CIHS}
Let $A\succ0_n$ be a Div-BV tensor over $\R^n$. Let $X_1,\ldots,X_m$ be points at which $A$ is, locally, Div-free and positively homogeneous of degree $1-n$. Then the functional inequality improves into
\begin{equation}\label{eq:FIDm}
\int_{\R^n}(\det A)^{\frac1{n-1}}dx+\sum_{j=1}^m{\rm Dm}(A;X_j)\le c_n\|{\rm Div}\,A\|_{\cal M}^{\frac n{n-1}},
\end{equation}
for some finite constant  $c_n$, independent from the tensor and the number $m$ of singularities.
\end{thm}

\subsection{Proof of Theorem \protect\ref{th:HS}}

Let $M$ be the mass-momentum tensor of the Hard Spheres configuration. It consists in a kinetic part $M_{\rm kin}$, supported by the trajectories of the centers of the particles, and a colliton part $M_{\rm col}$, described above. The interested reader will find all the details in the references \cite{Ser_HS,Ser_period}. We recall that $M=M_{\rm kin}+M_{\rm col}$ is Div-free.

As such, $M$ is useless in the context of Compensated Integrability, because it is essentially rank-one. Even at the nodes of $\cal G$, which $M$ does not charge, $M$ is only rank-two. Thus $(\det M)^{\frac1n}\equiv0$. This is why, following \cite{Ser_HS}, we add to $M$ correctors at each node. Their role is to make $M_1:=M+M_{\rm cor}$ full rank at the nodes, while remaining  locally Div-free and positively homogeneous. The price to pay is a small contribution to the Divergence, away from the nodes. We recall that the support of the correctors is arbitrarily small and, since the nodes do not approach from $\partial U$ by less than $r$, we may assume that the support of $M_{\rm cor}$ does not meet the boundary. Therefore $\gamma_\nu M_1=\gamma_\nu M$ and in particular $\gamma_\nu^\tau M_1\equiv0$.

So far, $M_1$ is Div-BV in $U=\R_+\times\Omega$ for every $T>0$. Choosing $T>0$, we define as usual a Div-BV tensor $A$ over $\R\times\Omega$ by $A=M_1$ if $t\in(0,T)$, and $A=0_n$ otherwise. The tangential part of the normal trace of $A$ still vanishes. According to (\ref{eq:estext}) and arguing as in Paragraph \ref{sss:evol}, there exists an extension $\A$ of $A$ to $\R^{1+d}$, with the property that 
$$\|{\rm Div}\,\A\|_{\cal M}\le_d\|{\rm Div}\,A\|_{\cal M}+K(\partial\Omega)\|A_{r}\|_{\cal M}.$$

We now apply Theorem \ref{th:CIHS} to $\A$ and the nodes $X_1,\ldots,X_m$ of $M$ in $(0,T)\times\Omega$, these being nodes of $M_1$ and therefore of $\A$ ($\A$ might have other nodes in $\R^n$ but we do not need to take them in account). From (\ref{eq:FIDm}) we obtain 
\begin{eqnarray*}
\sum_{j=1}^m{\rm Dm}(M_1;X_j) & \le & c_n\|{\rm Div}\,\A\|_{\cal M}^{1+\frac1d} \\
& \le_d & \left(mN+\sqrt{mNE_0\,}+\|{\rm Div}\,M_1\|_{\cal M}+K(\partial\Omega)\|(M_1)_{r}\|_{\cal M}\right)^{1+\frac1d},
\end{eqnarray*}
where the two first terms in the upper bound come, as usual, from the jumps of $A$ at initial and final time, which contribute to ${\rm Div}\,A$.

The rest of the proof follows closely that in \cite{Ser_period}. At each node $X=X_j$, the corrector is of the form $b_XS_X$ where $S_X$ is a normalized positive tensor and $b_X>0$ is a parameter to be chosen. We have
$${\rm Dm}(M_1;X)=2^{d-3}b_X^{1-\frac1d}\left|\binom m{mu_{\rm out}}\wedge\binom m{mu_{\rm in}}\right|^{\frac1d}\ge2^{d-3}b_X^{1-\frac1d}(m^2|u_{\rm out}-u_{\rm in}|)^{\frac1d}.$$
On another hand only the correctors contribute to $\|{\rm Div}\,M_1\|_{\cal M}$, by the amount $2(d-1)\sum b_X$. Eventually
$$|M_{br}|=m\sum_\pi|u|^2dt|_{\gamma(\pi)}+m\sum_{\rm coll.}|u_{\rm out}-u_{\rm in}|\,d\ell_J$$
and
$$|M_{tr}|=m\sum_\pi|u|dt|_{\gamma(\pi)}.$$
where two sums run over the $N$ particles, and the other runs over the collitons. 
This gives 
$$\|M_{br}\|_{\cal M}\le T\left(2E_0+\sqrt{2mNE_0\,}\right)+2rm\sum_{\rm coll.}^{(0,T)}|u_{\rm out}-u_{\rm in}|.$$
Eventually  the contribution of $M_{\rm cor}$ to the mass of $(M_1)_{br}$ can be taken arbitrarily small by shrinking the support of the correctors. There remains therefore the estimate
\begin{eqnarray*}
\sum_{\rm coll.}^{(0,T)}b_X^{1-\frac1d}(m^2|u_{\rm out}-u_{\rm in}|)^{\frac1d} & \le_d & \left(mN+\sqrt{mNE_0\,}+\sum_{\rm coll.}^{(0,T)} b_X \right. \\
& & +\left. K(\partial\Omega)(T\left(E_0+\sqrt{mNE_0\,}\right)+rm\sum_{\rm coll.}^{(0,T)}|u_{\rm out}-u_{\rm in}|)\right)^{1+\frac1d}.
\end{eqnarray*}

The next step is to set $b_X=\lambda\beta_X^{\frac d{d-1}}$, where $\beta_X$ is still to be chosen, and 
$$\lambda=\frac{mN+\sqrt{mNE_0\,}+K(\partial\Omega)(T\left(E_0+\sqrt{mNE_0\,}\right)+rm\sum_{\rm coll.}^{(0,T)}|u_{\rm out}-u_{\rm in}|)}{\sum_{\rm coll.}^{(0,T)} \beta_X^{\frac d{d-1}}}\,.$$
This gives
\begin{eqnarray*}
\sum_{\rm coll.}^{(0,T)}\beta_X(m^2|u_{\rm out}-u_{\rm in}|)^{\frac1d} & \le_d & \|\vec\beta\|_{\frac d{d-1}}\left(mN+\sqrt{mNE_0\,}\right. \\
& & \left.+K(\partial\Omega)(T\left(E_0+\sqrt{mNE_0\,}\right)+rm\sum_{\rm coll.}^{(0,T)}|u_{\rm out}-u_{\rm in}|)\right)^{\frac2d}.
\end{eqnarray*}
It is now the time to choose
$$\beta_X=(m^2|u_{\rm out}-u_{\rm in}|)^{\frac{d-1}d},$$
which yields
$$m^2\sum_{\rm coll.}^{(0,T)}|u_{\rm out}-u_{\rm in}|\le_d\left(mN+\sqrt{mNE_0\,}+K(\partial\Omega)\left(T\left(E_0+\sqrt{mNE_0\,}\right)+rm\sum_{\rm coll.}^{(0,T)}|u_{\rm out}-u_{\rm in}|\right)\right)^2.$$
The latter estimate is still unsatisfactory, because of its physical inhomogeneity. To overcome this flaw, we consider the initial configuration in which the particles have the same location $y_j(0)$, but have velocities $\mu u_j(0)$, where $\mu>0$ is some parameter. The motion is exactly the same, up to a time rescaling $t\mapsto \frac t\mu$\,. We apply our estimate to this new motion, on the time interval $(0,T/\mu)$. The kinetic energy of the new motion being $\mu^2E_0$, we obtain
\begin{eqnarray}
\nonumber
\mu m^2\sum_{\rm coll.}^{(0,T)}|u_{\rm out}-u_{\rm in}| & \le_d & \left(mN+\mu\sqrt{mNE_0\,} \right. \\
\label{eq:HSinH}
 & & \left.+ K(\partial\Omega)(T\left(\mu E_0+\sqrt{mNE_0\,}\right)+\mu rm\sum_{\rm coll.}^{(0,T)}|u_{\rm out}-u_{\rm in}|)\right)^2.
\end{eqnarray}
Chosing 
$$\frac1\mu=\frac1{mN+K(\partial\Omega)T\sqrt{mNE_0\,}}\,\left(\sqrt{mNE_0\,}+K(\partial\Omega)(E_0T+rm\sum_{\rm coll.}^{(0,T)}|u_{\rm out}-u_{\rm in}|)\right)$$
in (\ref{eq:HSinH}), we obtain a now homogeneous estimate
\begin{eqnarray}
\nonumber
m^2\sum_{\rm coll.}^{(0,T)}|u_{\rm out}-u_{\rm in}| & \le &  c_d\left(mN+K(\partial\Omega)T\sqrt{mNE_0\,}\right) \times \\
\label{eq:nowhom}
& & \left(\sqrt{mNE_0\,}+K(\partial\Omega)(E_0T+rm\sum_{\rm coll.}^{(0,T)}|u_{\rm out}-u_{\rm in}|)\right).
\end{eqnarray}

Let us now denote $\kappa_d=1/(3c_d)$,
where $c_d$ is the constant appearing in (\ref{eq:nowhom}). If 
$$NrK(\partial\Omega)<\kappa_d \qquad\hbox{and}\qquad T\le T_f:=\frac1{K(\partial\Omega)\bar u}\,,$$ 
then this inequality implies
$$\sum_{\rm coll.}^{(0,T)}|u_{\rm out}-u_{\rm in}|\le 3c_d^2N^2\bar u.$$
Since we may replace the initial time by an arbitrary $t\ge0$, this ends the proof of Theorem \ref{th:HS}.

\bigskip

\appendix

\section{Normal trace}\label{ap:trace}

Let ${\rm divBV}(U)$ denote the Banach space of vector fields $\vec q$ whose coordinates $q_i$, as well as the divergence, are finite measures over $U$. Notice the use of a capital letter in ${\rm DivBV}(U)$ when speaking of tensors, and of lower case when speaking of vector fields.

We prove that the normal trace operator 
\begin{eqnarray*}
{\rm divBV}(U) & \longrightarrow & (C^1(\partial U))' \\
\vec q & \longmapsto & \gamma_\nu \vec q
\end{eqnarray*}
is well-defined by the duality formula
\begin{equation}\label{eq:dual}
\forall Y\in C^1(\bar U),\qquad \langle \vec q,\nabla Y\rangle_U+\langle {\rm div}\,\vec q,Y\rangle_U=\langle \gamma_\nu \vec q,Y\rangle_{\partial U}.
\end{equation}
Since the map 
\begin{eqnarray*}
C^1(\bar U) & \to & C^1(\partial U) \\
Y & \mapsto & Y|_{\partial U}
\end{eqnarray*}
is onto (for a domain with $C^1$-boundary) and admits a bounded right inverse, it suffices to show that  the left-hand side in (\ref{eq:dual}) depends only upon the restriction $Y|_{\partial U}$, but not upon the normal derivative $(\vec\nu\cdot\nabla)Y$. Equivalently, that it vanishes whenever $Y|_{\partial U}\equiv0$. This is already true (by definition) when $Y$ is compactly supported. When $Y$ simply vanishes at the boundary, we approximate it by $Y_\epsilon=\phi_\epsilon Y$, where $\phi_\epsilon\in{\cal D}(U)$ satisfies $\phi_\epsilon\equiv1$ for $d(x)>\epsilon$, and $|\nabla\phi_\epsilon|\le\frac2\epsilon$ everywhere. Then 
\begin{eqnarray*}
|\langle \vec q,\nabla Y\rangle_U+\langle {\rm div}\,\vec q,Y\rangle_U| & = & |\langle \vec q,\nabla (Y-Y_\epsilon)\rangle_U+\langle {\rm div}\,\vec q,Y-Y_\epsilon\rangle_U| \\ 
& = & |\langle \vec q,(1-\phi_\epsilon)\nabla Y\rangle_U+\langle {\rm div}\,\vec q,(1-\phi_\epsilon)Y\rangle_U-\langle \vec q,Y\nabla \phi_\epsilon\rangle_U| \\
& \le & 3\|\nabla Y\|_\infty\|\vec q|_{V_\epsilon}\|_{\cal M}+\|Y\|_\infty\|({\rm div}\,\vec q)|_{V_\epsilon}\|_{\cal M},
\end{eqnarray*}
where the quantity $|Y|\cdot|\nabla\phi_\epsilon|$ is bounded by $2\|\nabla Y\|_\infty$ in $V_\epsilon$, the support of $\nabla\phi_\epsilon$.
Since both terms in the right-hand side tend to $0$ as $\epsilon\to0$, we obtain the announced property
$$\left(Y\in C^1(\bar U),\,Y|_{\partial U}\equiv0\right)\Longrightarrow\left(\langle \vec q,\nabla Y\rangle_U+\langle {\rm div}\,\vec q,Y\rangle_U=0\right).$$

\bigskip

If an additional function $f\in C^1(\bar U)$ is given, then 
\begin{eqnarray*}
\langle \gamma_\nu(f\vec q),Y\rangle & = & \langle f\vec q,\nabla Y\rangle+\langle {\rm div}\,(f\vec q),Y\rangle =  \langle \vec q,f\nabla Y\rangle+\langle f\,{\rm div}\,\vec q+\vec q\cdot\nabla f,Y\rangle \\ 
& = &   \langle \vec q,f\nabla Y\rangle+\langle {\rm div}\,\vec q,fY\rangle+\langle \vec q,Y\nabla f \rangle \\ 
& = & \langle \vec q,\nabla(f Y)\rangle+\langle {\rm div}\,\vec q,fY\rangle=\langle \gamma_\nu\vec q,fY\rangle=\langle f\gamma_\nu\vec q,Y\rangle.
\end{eqnarray*}
We conclude
\begin{equation}
\label{eq:Cungam}
\forall\,\vec q\in{\rm divBV}(U),\,\forall \,f\in C^1(\bar U),\qquad \gamma_\nu(f\vec q)=f\gamma_\nu\vec q.
\end{equation}

Eventually, the definition of $\gamma_\nu$ and the property (\ref{eq:Cungam}) extend to DivBV-tensors, since each rows are divBV-vector fields.

\section{The reflection across the boundary}\label{ap:refl}

Recall that if $\Sigma$ is a closed $C^3$-hypersurface in $\R^n$, then the unit normal vector 
\begin{eqnarray*}
\Sigma & \to & \R^n \\
a & \mapsto & \vec\nu
\end{eqnarray*}
is twice differentiable. Its differential at $a\in\Sigma$ maps the tangent space $T_a\Sigma$ into itself, defining a self-adjoint operator called the {\em curvature tensor}.

\bigskip

We suppose now that $\Sigma=\partial U$ is the $C^3$-boundary of the open domain $U$. Let $h:U\to\R_+$ denote the distance to $\Sigma$, a $C^2$-function in the corona $\omega$ defined by $h(x)<L=L(\Sigma)$, where we have $|\nabla h|\equiv1$. If $a\in\Sigma$ and $\mu\in(0,L)$, the operator
$$M_\mu:={\rm Id}_{T_a\Sigma}-\mu{\rm d}\vec\nu_a\in{\bf End}(T_a\Sigma)$$
is self-adjoint and invertible, with 
$$|M_\mu^{-1}|_{\rm op}\le\frac1{1-\mu\rho({\rm d}\vec\nu_a)}\,.$$
The spectral radius above is the largest absolute principal curvature of $\Sigma$ at $a$.

The projection $\pi:\omega\to\Sigma$ is given by $\pi(x)=x-h(x)\nabla h(x)$ and the reflection is 
$$\sigma(x)=x-2h(x)\nabla h(x).$$
Notice that $\vec\nu_{\pi(x)}=-\nabla h(x)$.
The symmetry of the differential 
$${\rm d}\,\sigma=I_n-2h{\rm D}^2h-2\nabla h\otimes\nabla h$$
reflects the fact that $\sigma$ is the gradient of a function $x\mapsto\frac12|x|^2-h(x)^2$.
\begin{lemma}
\label{lm:hddh}
Let $x\in\omega$ be given, and $a=\pi(x)$. Then 
$$h\,|{\rm D}^2h|_{\rm bil}=|h{\rm d}\vec\nu_a\left({\rm Id}_{T_a\Sigma}-h{\rm d}\vec\nu_a\right)^{-1}|_{\rm op}.$$
In particular $h\,|{\rm D}^2h|_{\rm bil}\le1$ and $|{\rm d}\,\sigma|_{\rm op}\le3$ in the sub-corona $\omega_{1/2}$ defined by $d(x;\Sigma)<L/2$. 
\end{lemma}

\bigskip

\bepr

We start from the identity $h(a-\mu\vec\nu_a)=\mu$. Differentiation gives
$$\nabla h\cdot(M_\mu\tau-\lambda\vec\nu_a)\equiv\lambda,$$
that is $\nabla h(a-\mu\vec\nu_a)=-\vec\nu_a$.

Differenciating once more in $\mu$, we find ${\rm D}^2h_x\vec\nu_a\equiv0$ where $a=\pi(x)$,  so that it suffices to evaluate ${\rm D}^2h_x(\tau,\tau)$ for vectors $\tau\in T_a\Sigma$. To do so, we differentiate in $a$ instead, to obtain
${\rm D}^2h\,M_\mu\tau= -{\rm d}\vec\nu_a\tau$.
Setting $\mu=h$, we obtain $h{\rm D}^2h\tau'=-K\tau'$ where
$$K=h{\rm d}\vec\nu_a\left({\rm Id}_{T_a\Sigma}-h{\rm d}\vec\nu_a\right)^{-1}.$$

\enpr

\bigskip

We shall restrict to the domain $\omega_{1/2}$ from now on, where we have also $|M_h^{-1}|_{\rm op}\le2$.
\begin{lemma}
\label{lm:ddh}
For $x\in\omega_{1/2}$ we have
$$|{\rm D}^2h_x|_{\rm bil}\le\frac2L\,.$$
\end{lemma}

\bigskip

\bepr

We have seen above ${\rm D}^2h_x\,\vec\nu_a=0$ and
$${\rm D}^2h_x\;\tau'=-{\rm d}\vec\nu_aM_h^{-1}\tau'.$$
Therefore
$$|{\rm D}^2h_x|_{\rm bil}\le|{\rm d}\vec\nu_a|_{\rm op}|M_h^{-1}|_{\rm op}\le\frac1L\,|M_h^{-1}|_{\rm op}\le\frac2L\,.$$

\enpr

\bigskip

As for the third differential we start on the one hand from ${\rm D}^2h(a-\mu\vec\nu_a)\,\vec\nu_a=0$. Differentiatiing in $\mu$, we have ${\rm D}^3h_x(\vec\nu_a,\vec\nu_a,\cdot)=0$. Differentiating instead in $a$, we obtain
$${\rm D}^3h_x(M_h\tau,\vec\nu_a,\cdot)-{\rm D}^2h_x({\rm d}\vec\nu_a\cdot\tau,\cdot)=0,$$
that is
$${\rm D}^3h_x(\tau',\vec\nu_a,\cdot)={\rm D}^2h_x({\rm d}\vec\nu_aM_h^{-1}\tau',\cdot).$$
Differentiating instead ${\rm D}^2h_xM_\mu\tau= -{\rm d}\vec\nu_a\tau$, we have
$${\rm D}^3h_x(M_h\tau,M_h\tau,\cdot)= (h{\rm D}^2h_x-1){\rm D}^2\vec\nu_a(\tau,\tau),$$
whence
$${\rm D}^3h_x(\tau',\tau',\cdot)= (h{\rm D}^2h_x-1){\rm D}^2\vec\nu_a(M_h^{-1}\tau',M_h^{-1}\tau').$$
Assembling the results above, we find
\begin{lemma}
For $x\in\omega_{1/2}$, we have
$$h(x)|{\rm D}^3h_x|_{\rm tri}\le\max\left(\frac2L\,,8L|{\rm D}^2\vec\nu_a|_{\rm bil}\right).$$
\end{lemma}

\bigskip

Since
$${\rm D}^2\sigma=-2\left(h{\rm D}^3h+{\rm D}^2h\otimes_{\rm sym}\nabla h\right),$$
we conclude therefore
\begin{cor}\label{c:ddsig}
For $x\in\omega_{1/2}$, we have
$$|{\rm D}^2\sigma|_{\rm bil}\le16\max\left(\frac1L\,,L|{\rm D}^2\vec\nu_a|_{\rm bil}\right).$$
\end{cor}


\begin{thebibliography}{00}




\bibitem{ADHR} A. Arroyo-Rabasa, G. De Philippis, J. Hirsch, F. Rindler. Dimensional estimates and rectifiability for measures satisfying linear PDE constraints. {\em Geom. Func. Anal.}, {\bf29} (2019), pp 639--658.


\bibitem{Bab} V. M. Babich. On the extension of functions (Russian). {\em Uspekhi Mat. Nauk}, {\bf8}, 2 (1953), pp 111--113





\bibitem{Pog} A. V. Pogorelov. {\em The Minkowski multidimensional problem}. Scripta Series in Mathematics. V. H. Winston \& Sons, Washington, D.C.; Halsted Press (John Wiley \& Sons), New York--Toronto--London (1978).

       
\bibitem{Ser_DPT} D. Serre. Divergence-free positive symmetric tensors  and fluid dynamics. {\em Annales de l'Institut Henri Poincar\'e (analyse non lin\'eaire)}, {\bf35} (2018), pp 1209--1234.

\bibitem{Ser_JMPA} D. Serre. Compensated integrability. Applications to the Vlasov--Poisson equation and other models of mathematical physics. {\em J. Math. Pures \& Appl.}, {\bf127} (2019), pp 67--88.

\bibitem{Ser_HS} D. Serre. Hard spheres dynamics: Weak {\em vs} strong collisions. {\em Arch. Rat. Mech. Anal.}, {\bf240} (2021), pp 243--264.
 
\bibitem{Ser_period} D. Serre. Compensated integrability on tori;  {\em a priori} estimate for space-periodic gas flows. {\em Comptes Rendus, Math\'ematiques}, to appear.
 
\bibitem{Tem} R. Temam. {\em Navier-Stokes equations: theory and numerical analysis.} North-Holland, Amsterdam (1977).
 


\end{thebibliography}
\end{document}